\documentclass[11pt]{amsart}
\usepackage{amsmath,amssymb}

\begin{document}

\title[critical condition]{  A critical regularity condition on the angular velocity of 
 axially symmetric  Navier-Stokes equations}
\thanks{}
\thanks{AMS Subject Classifications:  35Q30 and 35B07.}
\date{2015 May 1st }

\vspace{ -1\baselineskip}

\author{  Qi S. Zhang}

\address{Department of Mathematics,  University of California,
Riverside, CA 92521}

\numberwithin{equation}{section}
\allowdisplaybreaks

\newtheorem{theorem}{Theorem}
\newtheorem{proposition}{Proposition}
\newtheorem{lemma}{Lemma}
\newtheorem{definition}{Definition}
\newtheorem{corollary}{Corollary}
\newtheorem{remark}{Remark}[section]
\newtheorem{acknowledgement}{Acknowledgement}
\numberwithin{theorem}{section}
\numberwithin{lemma}{section}
\numberwithin{proposition}{section}
\numberwithin{equation}{section}
\numberwithin{definition}{section}
\def\R{\bf R}
\def\eps{\varepsilon}
\def\b{\bar}
\def\al{\aligned}
\def\eal{\endaligned}
\def\be{\begin{equation}}
\def\ee{\end{equation}}
\def\lab{\label}
\def\nn{\nonumber}
\def\~{\widetilde}
\def\->{\overrightarrow}
\def\bar{\overline}
\def\pd{\partial}
\def\nab{\nabla}
\def\lam{\lambda}

\newcommand{\reals}{\mathbb{R}}
\newcommand{\Div}{\textrm{div }}
\newcommand{\Curl}{\textrm{curl }}
\newcommand{\supp}{\textrm{supp }}
\newcommand{\V}{\widetilde{v}}
\newcommand{\Omegabar}{\overline{\Omega}}
\newcommand{\e}{\epsilon}
\newcommand{\om}{\omega}
\newcommand{\Om}{\Omega}

 \smallskip

\begin{abstract}
Let $v$ be the velocity of Leray-Hopf solutions to the axially
symmetric three-dimensional Navier-Stokes equations. 
It is shown that $v$ is regular if the angular velocity $v_\theta$ satisfies an integral 
condition which is critical under the standard scaling. This condition allows functions 
satisfying 
\[
|v_\theta(x, t)| \le \frac{C}{r |\ln r|^{2+\e}}, \quad r<1/2,
\] where $r$ is the distance 
from $x$ to the axis, $C$ and $\e$ are any positive constants.   

Comparing with the critical a priori bound 
\[
|v_\theta(x, t)| \le \frac{C}{r},  \qquad 0< r \le 1/2,
\]our condition is off by the log factor   $|\ln r|^{2+\e}$ at worst. 
This is inspired by the recent interesting paper \cite{CFZ:1} where H. Chen, D. Y. Fang and T. Zhang
 establish, among 
other things, an almost critical regularity condition on the angular velocity.  Previous regularity conditions
are off by a factor $r^{-1}$.

The proof is based on the new observation that, when viewed differently,  all the vortex stretching terms in
 the 3 dimensional axially symmetric Navier-Stokes equations are  
critical instead of supercritical as commonly believed.
\end{abstract}
\maketitle

\maketitle
\tableofcontents

\section{Introduction}

In
rectangular coordinates,
the incompressible Navier-Stokes equations are
\be
\lab{3dNS}
\Delta v-(v\cdot\nabla)v-\nabla p-\partial_t v=0,\ \Div v=0,
\ee
where
$v=(v_1(x,t),v_2(x,t),v_3(x,t)):\reals^3\times[0,T]\rightarrow\reals^3$
is the velocity field and
$p=p(x,t):\reals^3\times[0,T]\rightarrow\reals$ is the pressure.
In cylindrical coordinates $r,\theta, x_3$ with
$(x_1,x_2,x_3)=(r\cos\theta,r\sin\theta, x_3)$,  axially symmetric
solutions are of the form
\begin{align*}
v(x,t)=v_r(r,x_3,t)\->{e_r}+v_{\theta}(r,x_3,t)
\->{e_{\theta}}+v_3(r,x_3,t)\->{e_3}.
\end{align*}
The components $v_r,v_{\theta},v_3$
 are all independent of the angle of rotation $\theta$.
 Here $\->{e_r},\->{e_{\theta}},\->{e_3}$ are the basis vectors for
 $\reals^3$ given by
\begin{align*}
\->{e_r}= \Big (\frac{x_1}{r},\frac{x_2}{r},0 \Big ),
\ \->{e_{\theta}}= \Big (\frac{-x_2}{r},\frac{x_1}{r},0 \Big ),
\ \->{e_3}=(0,0,1).
\end{align*}

It is known   (see
\cite{CSTY1} for example) that $v_r$, $v_3$ and $v_\theta$ satisfy the equations
\begin{align}
\lab{eqasns}
\begin{cases}
   \big (\Delta-\frac{1}{r^2} \big )
v_r-(b\cdot\nabla)v_r+\frac{v_{\theta}^2}{r}-\partial_r
p-\partial_t v_r=0,\\
   \big   (\Delta-\frac{1}{r^2}  \big
)v_{\theta}-(b\cdot\nabla)v_{\theta}-\frac{v_{\theta}v_r}{r}-\partial_t v_{\theta}=0,\\
 \Delta v_3-(b\cdot\nabla)v_3-\partial_3 p-\partial_t v_3=0,\\
 \frac{1}{r}\partial_r (rv_r) +\partial_3
v_3=0,
\end{cases}
\end{align}
 where $b(x,t) =(v_r,0,v_3)$ and
the last equation is the divergence-free
 condition. Here, $\Delta$ is the cylindrical,  scalar Laplacian
and $\nabla$ is the cylindrical gradient field:
 \begin{align*}
 \Delta =\partial^2_r+
\frac{1}{r} \partial_r+\frac{1}{r^2}\partial
^2_\theta+\partial^2_3,\ \
 \nabla =  \Big  (\partial_r, \frac{1}{r} \partial_\theta, \pd_3  \Big  ).
\end{align*}
Observe that the equation for $v_{\theta}$ does not depend on the
pressure. Let $\Gamma=rv_{\theta}$, then
\begin{equation}\label{Gamma/vtheta}
\Delta \Gamma -(b\cdot\nabla)\Gamma-\frac{2}{r} \pd_r
\Gamma-\pd_t \Gamma=0,\ \Div b=0.
\end{equation}
The vorticity  $\omega=\Curl v$ for axially
symmetric solutions
 \begin{align*}
 \omega(x,t)=\omega_r\->{e_r}+\omega_{\theta}
\->{e_{\theta}}+\omega_3\->{e_3}
\end{align*} is given by
\begin{equation}
\label{curlformulas} \omega_r =-\pd_3 v_{\theta}, \ \omega_ {\theta}
 =\pd_3 v_r-\pd_r v_3, \ \omega_3=\pd_r v_
{\theta}+\frac{v_{\theta}}{r}.
\end{equation}
The equations of vorticity $\omega= \Curl v$ in
cylindrical form are  (again, see \cite{CSTY1} for example):
 \begin{align}
\lab{eqvort}
\begin{cases}
  \big (\Delta-\frac{1}{r^2} \big
)\omega_r-(b\cdot\nabla)\omega_r+\omega_r
\pd_r v_r +\omega_3\pd_3 v_r -\pd_t
\omega_r
=0,\\
   \big  (\Delta-\frac{1}{r^2}  \big
)\omega_{\theta}-(b\cdot\nabla)\omega_{\theta}+2\frac{v_{\theta}}
{r}\pd_3 v_{\theta}+\omega_{\theta}\frac{v_r}{r}-\pd_t
\omega_{\theta}=0,\\
 \Delta\omega_3-(b\cdot\nabla)\omega_3+\omega_3\pd_3
v_3+\omega_{r}\pd_r v_3 -\pd_t \omega_3=0.
\end{cases}
\end{align}

Although the axially symmetric Navier-Stokes equations is a special case of the full $3$ 
dimensional one, our level of understanding had been roughly the same, with essential 
difficulty unresolved. One quick explanation of the difficulty goes as follows. Viewing 
(\ref{3dNS}) as a reaction diffusion equation. The standard theory for regularity requires 
the velocity to be bounded in suitable function space whose norm is invariant under
standard scaling, such as $L^{p, q}$ with $\frac{3}{p}+\frac{2}{q}=1$.
 However the only general a priori bound available is the energy estimate, which 
scales as $-1/2$. So there is a positive gap between the two which makes the equations supercritical.

Equation (\ref{eqasns}) has been studied by many authors in recent years. The following is 
a list which is far from complete. If the swirl $v_\theta=0$, 
then  long time ago,  O. A. Ladyzhenskaya
\cite{L}, M. R. Uchoviskii and B. I. Yudovich \cite{UY}),  proved that finite energy solutions
to (\ref{eqasns}) are smooth for all time. See also the paper by S.
Leonardi, J. Malek, J. Necas,  and  M. Pokorny \cite{LMNP}).

In the presence of swirl, it is not known in general if finite energy solutions
blow up in finite time.   However a lower bound for the possible blow up rate is known by the recent results of
C.-C. Chen, R. M. Strain, T.-P.Tsai,  and  H.-T. Yau in \cite{CSTY1},
\cite{CSTY2},  G. Koch, N. Nadirashvili, G. Seregin,  and  V.
Sverak in \cite{KNSS}. See also the work by G. Seregin  and
V. Sverak \cite{SS} for a localized version. These authors prove that if
$
|v(x, t)| \le
\frac{C}{r},
$ then solutions are smooth for all time. Here $C$ is any positive constant. 
Their result can be rephrased as: type I solutions are regular.
See also the papers  \cite{LZ11}, \cite{LZ11-2} on further results in this direction.
 J. Neustupa  and
M. Pokorny \cite{NP} proved that the regularity of one component
(either $v_r$ or $v_{\theta}$) implies regularity of the other
components of the solution. See more refined results in \cite{NP:2} and 
the work of Ping Zhang and Ting Zhang \cite{ZZ}. Also proving regularity is the work of
Q. Jiu  and  Z. Xin \cite{JX} under an assumption of sufficiently small
zero-dimension scaled norms.  D. Chae  and  J. Lee \cite{CL}  also proved
regularity results assuming finiteness of another certain
zero-dimensional integral.  G. Tian  and  Z.
Xin \cite{TX}  constructed a family of singular axially
symmetric solutions with singular initial data.
T. Hou  and  C. Li \cite{HL} found a special class of global smooth
solutions. See also a recent extension: T. Hou, Z. Lei  and  C. Li
\cite{HLL}.

 Define
\[
J= \frac{\omega_r}{r},  \quad \Omega=\frac{\omega_{\theta}}{r}.
\]Then the triple $J, \Om, \om_3$ 
 satisfy the system
\begin{equation}
\label{eqjoo}
\begin{cases}
\Delta J  -(b\cdot\nabla) J +\frac{2}{r}\pd_r J +
 (\om_r \pd_r + w_3 \pd_3) \frac{v_r}{r} - \pd_t J
=0,\\
\Delta \Omega -(b\cdot\nabla)\Omega+\frac{2}{r}\pd_r
\Omega - \frac{2v_{\theta}}{r} J -\pd_t \Omega=0,\\
\Delta w_3-(b\cdot\nabla) w_3+ w_{r}\pd_r v_3 + w_3 \pd_3
v_3 -\pd_t w_3=0. 
\end{cases}
\end{equation}  Here, in the second equation, we used the identity $r J = w_r = - \pd_3 v_\theta$. 

A great observation by Hui Chen, Daoyuan Fang and Ting Zhang in \cite{CFZ:1} is that 
the first two equations in (\ref{eqjoo}) form a critical system under the 
standard scaling.  Using this and a "magic formula" relating $\nabla (v_r/r)$ with $w_\theta/r$ by 
Changxing Miao and Xiaoxin Zheng \cite{MZ}, they obtained, among other things, an almost critical regularity 
condition on $v_\theta$.  For example it is proven that if $|v_\theta(x, t)| \le C/r^{2-\e}$ with
$\e>0$, then solutions are regular.

In this paper we observe further that, all three equations are critical
 when viewed in a suitable way. Therefore the vorticity equation of
 3 dimensional axially symmetric Navier-Stokes equations are  
critical instead of supercritical as commonly believed. 
This, together with a localization method in \cite{Z}, allow us to prove 
Theorem \ref{mainthm} below, which provides a localized critical regularity condition on $v_\theta$. 
It is tantalizing that our condition differs with the critical a priori bound (\cite{CL} or \cite{NP})
\[
|v_\theta(x, t)| \le \frac{C}{r},  \qquad 0< r \le 1/2,
\] by the log factor   $|\ln r|^{2+\e}$ at worst. See the remarks below. 

Now we introduce the function class where $v_\theta$ lives. It is defined in an integral way 
which is usually called the form boundedness condition, which is more general than the 
corresponding $L^{p, q}$
condition.

\begin{definition}  
\lab{defcrit}
We say the angular velocity $v_\theta$ is in the $\lam_1$ critical class if 
there is a positive number $a <1$ and another positive number $\lambda_2$ such
that the inequality  
\[
\int^t_0\int \left(\frac{|v_\theta|}{r} + v^2_\theta \right) \psi^2 dyds 
\le \lam_1
\int^t_0\int |\nab \psi|^2 dyds + \frac{\lam_2}{a^2}
\int^t_0\int \psi^2 dyds 
\]holds for all $t \ge 0$ and for all smooth $\psi=\psi(y, s)$, $s \in [0, t]$,  satisfying the conditions
(1) $\psi$ is axially symmetric in $y$; (2)  $\psi(\cdot, s)$ is supported in the 
cylinder $D_{a, l}=\{ (r, \theta, x_3) \, | \,  0 \le r <a,  \, -l < x_3 < l, 
0 \le \theta < 2 \pi  \}$  for some $l \ge a$.
\end{definition}

\begin{remark}
\rm Clearly the class is scaling invariant. A function
$v_\theta$ is  the $\lam_1$ critical class for all $\lam_1>0$ if it satisfies 
$|v_\theta(x, t)| \le \frac{C}{r |\ln r|^{2+\e}}, \quad r<1/2$. Here $C>0$, $\e>0$ are 
arbitrary positive constant.   This claim will be proven at the end of the paper.
One may also take $\e=0$ but replace $r$ by $r/a$ and $C$ by a small constant in the bound, by
virtue of the 2 dimensional Hardy's inequality.
\end{remark}

Here is the main result of the paper.

\begin{theorem} \label{mainthm}
 Let  $v$ be a Leray-Hopf
axially symmetric solution of the three-dimensional Navier-Stokes
equations in $\reals^3\times(0, \infty)$ with initial data
$v_0=v(\cdot, 0)\in L^2(\reals^3)$. Assume further
$rv_{0,{\theta}}\in L^{\infty}(\reals^3)$.  

There exists a positive number $\lam_1$.
Suppose $v_\theta$ is in the $\lam_1$ critical class. 
Then $v$ is smooth for all time.
\end{theorem}

\begin{remark}
\rm The size of $\lam_1$ is estimated in (\ref{lam1}). It is an absolute 
constant depending on the $L^2$ norm of the Riesz operators. There is no size restriction 
on $\lam_2$. Also the $a^2$ in the definition can be replaced by any positive continuous function 
of $a$.  But this may break the scaling invariance.
\end{remark}

The theorem will be proven in the next section.
The following are some notations to be frequently used. We use $x=(x_1,x_2,x_3)$ to denote
a point in $\reals^3$ for rectangular coordinates, and in the
cylindrical system we use $r=\sqrt{x_1^2+x_2^2}$,
$\theta=\tan^{-1}\frac{x_2}{x_1}$.
We will use $S(v_0, ...), C(v_0, ...)$ to denote positive constants which depend on the 
initial velocity $v_0$ etc. Also $C$ denotes absolute constant which may change value. 

Let us explain why the vortex stretching terms in (\ref{eqjoo}) are critical.  
For example the term $w_3 \pd_3 v_3$ where $\pd_3 v_3$ being viewed as a potential of the 
unknown function $w_3$ is certainly supercritical. However, we view $w_3 = \pd_r v_\theta +
\frac{v_\theta}{r}$ as the potential and $\pd_3 v_3$ as the unknown.  Since it is known that 
$|v_\theta| \le C/r$, we see that $w_3$ now scales as $-2$ power of the distance. 
This scaling shows $w_3$
 is a critical potential function.  The unknown function
$\pd_3 v_3$ scales the same way as the vorticity $w$. By exploiting the integral relations 
between $v$ and $w$, we can convert $\pd_3 v_3$ into $w_r, w_3, w_\theta$.  This, combined 
with the observation \cite{CFZ:1} about the first two equations in (\ref{eqjoo}),  imply that all the 
vortex stretching terms are critical.  Next we carry a local energy estimate for 
$(J, \Om, w_z)$ via equations (\ref{eqjoo}).  Once we know the potential terms 
are critical,  the drift terms can be treated by an old  small trick in \cite{Z},
the proof thus goes through.

\section{Proof of the theorem}

The proof is divided into several steps. We may assume that $v$ is smooth up to  a given time 
$t$.
\medskip

{\it Step 1. Choose suitable test functions for equations (\ref{eqjoo}).}

It is well known that singularity can possibly appear only on a finite segment of 
the $x_3$ axis  (\cite{CKN:1} for suitable solutions and \cite{BZ:1} for general ones). So by picking any positive number $a \le 1$ and another 
positive number $l>a$, which may depend on the initial velocity $v_0$, we can ensure that
$v$ is regular outside of the domain 
$D_1= \{ (r, \theta, x_3) \, | \,  0 \le r <a/2,  \, -l/2 < x_3 < l/2, 0 \le \theta < 2 \pi  \}$ for all time.
Let $\phi=\phi(r, x_3)$ be a axially symmetric cut off function in 
$D_2=\{ (r, \theta, x_3) \, | \,  0 \le r <a,  \, -l < x_3 < l, 
0 \le \theta < 2 \pi  \}$ such that $\phi=1$ on 
$D_3= \{ (r, \theta, x_3) \, | \,  0 \le r <2a/3,  \, -2l/3 < x_3 < 2l/3, 
0 \le \theta < 2 \pi  \}$ and $\phi=0$ on $D_2^c$ and also 
$\frac{|\nabla \phi|}{\phi^{1/2}} \le C/a$, $|\nabla^2 \phi| \le 
C/a^2$. 

Use $J \phi^2$, $\Om \phi^2$ and $w_3 \phi^2$ as test functions in equations 1, 2 and 3 in (\ref{eqjoo})
respectively. After integration on the region $D_2 \times [0, t]$ for $t>0$ we find that
\be
\lab{eqJ}
\al
L_1 &\equiv -\int^t_0 \int \Delta J \,  J \phi^2 dy ds - \int^t_0 \int \frac{2}{r} \pd_r J \, J \phi^2 dyds 
+ \int^t_0 \int \pd_t J \, J \phi^2 dyds \\
&= - \int^t_0 \int b \nabla J \, J \phi^2 dyds 
+ \int^t_0 \int (w_r \pd_r \frac{v_r}{r} + w_3 \pd_3 \frac{v_r}{r})
J \phi^2 dyds\\
&\equiv R_1 + T_1.
\eal
\ee

\be
\lab{eqOm}
\al
L_2 &\equiv -\int^t_0 \int \Delta \Om \,   \Om \phi^2 dy ds - \int^t_0 \int \frac{2}{r} \pd_r  \Om \,  \Om \phi^2 dyds 
+ \int^t_0 \int \pd_t  \Om \,  \Om \phi^2 dyds \\
&= - \int^t_0 \int b \nabla  \Om \,  \Om \phi^2 dyds 
- \int^t_0 \int \frac{ 2 v_\theta}{r} J \Om \phi^2 dyds\\
&\equiv R_2 + T_2.
\eal
\ee

\be
\lab{eqW3}
\al
L_3 &\equiv -\int^t_0 \int \Delta w_3  \,  w_3  \phi^2 dy ds 
+ \int^t_0 \int \pd_t w_3  \, w_3  \phi^2 dyds \\
&= - \int^t_0 \int b \nabla w_3  \, w_3  \phi^2 dyds +
 \int^t_0 \int (w_3 \pd_3 v_3 + w_r \pd_r v_3)
 w_3  \phi^2 dyds\\
&\equiv R_3 + T_3.
\eal
\ee

The left hand side of the three equalities $L_1$, $L_2$ and $L_3$ can be treated by routine integration by parts
which shows:
\[
\al
L_1 &= \int^t_0 \int | \nab J |^2 \phi^2 dyds + \int^t_0 \int J^2(0, y_3, t) \phi^2 dy_3 d r d t + 
\frac{1}{2} \int J^2 \phi^2 dy \bigg |^t_0 \\
& \qquad-  
\int^t_0 \int \nab J \, J  \nab \phi^2 dyds +  \int^t_0 \int  J^2 \frac{\pd_r \phi^2}{r} dyds.
\eal
\]Therefore
\[
\al
L_1 &\ge \frac{1}{2}  \int^t_0 \int | \nab J |^2 \phi^2 dyds  + 
\frac{1}{2} \int J^2 \phi^2 dy \bigg |^t_0 \\
& \qquad-  
2 \int^t_0 \int  J^2  |\nab \phi|^2 dyds +  \int^t_0 \int  J^2 \frac{\pd_r \phi^2}{r} dyds.
\eal
\]By our choice of the cut off function $\phi$, we know $v$ is regular in the supports of $\nab \phi$ and $\pd_r \phi$,
 which is bounded away from the singular set by a distance $a/6$. So there is a positive constant $S=S(v_0, a, l)$ such that
\be
\lab{l1}
L_1 \ge \frac{1}{2}  \int^t_0 \int | \nab J |^2 \phi^2 dyds  + 
\frac{1}{2} \int J^2 \phi^2 dy \bigg |^t_0- C t S(v_0, a, l).
\ee  Here we recall that $J$ and $\Om$ are all smooth functions if $v$ is smooth. 
Similarly
\be
\lab{l2}
L_2 \ge \frac{1}{2}  \int^t_0 \int | \nab \Om |^2 \phi^2 dyds  + 
\frac{1}{2} \int \Om^2 \phi^2 dy \bigg |^t_0- C t S(v_0, a, l),
\ee
\be
\lab{l3}
L_3 \ge \frac{1}{2}  \int^t_0 \int | \nab w_3 |^2 \phi^2 dyds  + 
\frac{1}{2} \int w_3^2 \phi^2 dy \bigg |^t_0- C t S(v_0, a, l).
\ee We remark that $S(v_0, a, l)$ may blow up when $a \to 0$. But we will make $a$ small and fixed.  

Substituting(\ref{l1}), (\ref{l2}) and (\ref{l3}) into (\ref{eqJ}), (\ref{eqOm}) and (\ref{eqW3})
respectively, we deduce
\be
\lab{jowz}
\al
&\int \left(J^2+ \Om^2 +w^2_3 \right) \phi^2 dy \bigg |^t_0
+   \int^t_0 \int \left( | \nab J|^2 +|\nab \Om|^2 + | \nab w_3 |^2 \right) \phi^2 dyds\\
&\le 2 (R_1+ R_2+R_3) + 2 (T_1+T_2+T_3) + C S(v_0, a, l).
\eal
\ee We are going to bound the right hand side in the next few steps.

\medskip

{\it Step 2.  bounds on $R_1+ R_2+R_3$, the drift terms. }

These terms are generated by $b= v_r \-> e_r + v_3 \-> e_3$ which is supercritical.
 However since these are given by divergence free drift terms, they can be bounded as done in \cite{Z}. We present a proof for completeness. 

Since $\Div b=0$, we have
\begin{align*} R_1&=-\int^t_0 \int b\cdot(\nabla J)(J
\phi^2)dyds\\
&= \int^t_0 \int b\cdot(\nabla\phi)\phi J^2dyds\\
&\leq    \Big  |\int   \Big  (
b\phi^{3/2}|J|^{3/2}   \Big  )  \Big (\frac{\nabla
\phi}{\phi^{1/2}}|J|^{1/2}  \Big )dyds    \Big  |.
\end{align*}
By H\"{o}lder's
inequality with exponents $\frac{4}{3}$ and 4,
\begin{align*}
R_1&\leq    \Big   (\int^t_0 \int |b|^{\frac{4}{3}}  \Big
(\phi^{3/2}|J|^{3/2}   \Big  )^{\frac{4}{3}}dyds
 \Big )^{\frac{3}{4}}   \Big  (\int^t_0\int  \Big
(\frac{|\nabla
\phi|}{\phi^{1/2}}|J|^{1/2}   \Big  )^4dyds  \Big
)^{\frac{1}{4}}.
\end{align*}
Using properties of the cutoff function we find:
   \begin{align*}
R_1\leq   \Big (\int^t_0 \int |b|^{\frac{4}{3}}(J\phi)^2dyds
 \Big  )^{\frac{3}{4}}\frac{C }{a}  \Big
(\int^t_0\int_{supp \, |\nab \phi|} J^2dyds  \Big )^{\frac{1}{4}}.
\end{align*}
 Next we fix
$\epsilon_1>0$ and we apply Young's inequality, with exponents
$\frac{4}{3}$ and $4$:
 \begin{align*}
 R_1&\leq  \Big (\frac{4}{3}\epsilon_1
 \Big  )^{\frac{3}{4}}  \Big(
\int^t_0 \int |b|^{\frac{4}{3}}(J \phi) ^2dyds  \Big
)^{\frac{3}{4}}\cdot  \Big  (\frac{4}{3}\epsilon_1  \Big
)^{-\frac{3}{4}}\frac{C }{a} 
\Big(\int^t_0 \int_{supp \, | \nab \phi |}  J ^2dyds  \Big )^{\frac{1}{4}}\\
&\leq\epsilon_1 \int^t_0\int |b|^{\frac{4}{3}}
(J \phi)^2dyds+\frac{C \epsilon_1^{-3}}
{a^4}\int^t_0 \int_{supp \, | \nab \phi |} J ^2dyds.
\end{align*}
 Thus,
\begin{equation}
|R_1|\leq
\epsilon_1c_0  \Vert b \Vert_{^{2, \infty}}^{4/3} 
\int^t_0 \int |\nabla(J\phi)|^2dyds +\frac{C \epsilon_1^{-3}}
{a^4}\int^t_0 \int_{supp \, | \nab \phi |} J ^2dyds.
\end{equation} 
  This
last inequality holds as a result of the standard energy estimate, H\"{o}lder's inequality
with exponents $\frac{3}{2}$ and $3$, and the $3$ dimensional Sobolev Inequality,
\begin{align*}
\int^t_0 \int |b|^{\frac{4}{3}}(J \phi)^2dyds&
\leq\int^t_0  \Big (\int |b|^2dy  \Big )^{\frac{2}{3}}
 \Big (\int (J \phi)^6dy   \Big )^ {\frac{1}{3}}ds\\
&\leq c_0  \Vert b \Vert_{^{2, \infty}}^{4/3}  \int^t_0 \int |\nabla(J \phi)|^2dyds.
\end{align*}
By choosing $\e_1$ suitably, we deduce
\be
\lab{R1<}
|R_1|\leq
 \frac{1}{8} \int^t_0 \int |\nabla J|^2 \phi^2 dyds +
C S(v_0, a, l),
\ee where we have used the fact that $v$ is regular in the support
of $\nab \phi $ for all time. In exactly the same manner, we 
find that
\be
\lab{R123<}
|R_1| + |R_2| + |R_3| \leq
 \frac{1}{8} \int^t_0 \int  \left( |\nabla J|^2  +
|\nabla \Om|^2 + |\nabla w_3|^2 \right) \phi^2 dyds +
C S(v_0, a, l),
\ee

\medskip

{\it Step 3.  bounds on $T_1$ and $T_2$.} 

In this step we follow the idea in [CFZ] with one modification, namely 
a localized version of a formula of Miao and Zheng which relates $\frac{v_r}{r}$ with 
$\frac{w_\theta}{r}$.  The rest of the step is divided into a few sub steps.

{\it step 3.1}

First we work on the easy one $T_2$ defined in (\ref{eqOm}).
\[
\al
T_2 &=  - \int^t_0 \int \frac{ 2 v_\theta}{r} J \Om \phi^2 dyds\\
& \le \int^t_0 \int \frac{|v_\theta|}{r} (J \phi)^2 dyds + 
           \int^t_0 \int \frac{|v_\theta|}{r} (\Om \phi)^2 dyds.
\eal
\]By our assumption on $v_\theta$, this implies
\[
T_2 \le \lambda_1 \int^t_0 \int  (|\nab (J \phi)|^2+ |\nab (\Om \phi)|^2) dyds + 
           \lambda_2 \int^t_0 \int [ (J \phi)^2 +  (\Om \phi)^2] dyds.
\]Let us write $\nab (J \phi) = \nab J \phi + J \nab \phi$. As mentioned earlier, $J$ 
is regular in the support of $\nab \phi$. Hence  
\be
\lab{T2<} 
T_2 \le 2 \lambda_1 \int^t_0 \int  (|\nab J|^2+ |\nab \Om |^2) \phi^2 dyds + 
           \lambda_2 \int^t_0 \int [ (J \phi)^2 +  (\Om \phi)^2] dyds + C t S(v_0, a, l).
\ee Here we also did the same argument for $\nab (\Om \phi)$.

{\it step 3.2}

Next we turn to $T_1$.  From (\ref{eqJ}),
\[
\frac{d T_1}{d t}  =  \int   (w_r \pd_r \frac{v_r}{r} + w_3 \pd_3 \frac{v_r}{r})
J \phi^2 dy
\]Using the relation $w_r = - \pd_3 v_\theta$, $w_3 = \frac{1}{r} \pd_r(r v_\theta)$
and integration by parts,
we see that
\[
\al
\frac{d T_1}{d t} &= - \int \pd_3 v_\theta \pd_r (\frac{v_r}{r}) J \phi^2 dy + 
\int  
\frac{1}{r} \pd_r(r v_\theta) \pd_3 (\frac{v_r}{r}) J \phi^2 dy\\
&=\int v_\theta \pd_3 \pd_r (\frac{v_r}{r}) J \phi^2 dy + 
\int v_\theta  \pd_r (\frac{v_r}{r}) \pd_3 (J \phi^2) dy\\
&\qquad - \int v_\theta \pd_r \pd_3 (\frac{v_r}{r}) J \phi^2 dy 
- \int v_\theta \pd_3 (\frac{v_r}{r})  \pd_r (J \phi^2) dy.
\eal
\]Notice that the first and third term on the right hand side of the last equality cancel.
Therefore,  we deduce
\[
\al
\frac{d T_1}{d t} 
&= 
\int v_\theta  \pd_r (\frac{v_r}{r}) (\pd_3 J) \phi^2 dy-
 \int v_\theta \pd_3 (\frac{v_r}{r})  (\pd_r J) \phi^2 dy\\
&\qquad +\int v_\theta  \pd_r (\frac{v_r}{r})  J \pd_r \phi^2 dy
- \int v_\theta \pd_3 (\frac{v_r}{r})  J \pd_r \phi^2 dy.
\eal
\]This implies, since the last two terms in the above identity are bounded, that
\[
\al
 T_1
&\le \frac{1}{8} \int^t_0 \int |\pd_3 J|^2 \phi^2 dy 
+ 2  \int^t_0 \int v_\theta^2 | \pd_r \frac{v_r}{r} |^2 \phi^2 dy \\
&\qquad +\frac{1}{8}  \int^t_0 \int |\pd_r J|^2 \phi^2 dy 
+ 2  \int^t_0 \int v_\theta^2 | \pd_3 \frac{v_r}{r} |^2 \phi^2 dy + C t S(v_0, a, l).
\eal
\]By our condition on $v_\theta$ again, we find that
\[
\al
 T_1
&\le \frac{1}{8} \int^t_0 \int |\nab J|^2 \phi^2 dy + C t S(v_0, a, l)
+ 2 \lam_1 \int^t_0 \int  | \nab ( \phi \pd_r \frac{v_r}{r}) |^2  dy
+ 2 \lam_2 \int^t_0 \int   ( \phi \pd_r \frac{v_r}{r})^2  dy\\
&\qquad + 2 \lam_1 \int^t_0 \int  | \nab ( \phi \pd_3 \frac{v_r}{r}) |^2  dy
+ 2 \lam_2 \int^t_0 \int   ( \phi \pd_3 \frac{v_r}{r})^2  dy.
\eal
\]This implies, after using again the fact that $v$ is smooth in the support 
of $\nab \phi$, that
\be
\lab{T1<4}
\al
 T_1
&\le \frac{1}{8} \int^t_0 \int |\nab J|^2 \phi^2 dy + C t S(v_0, a, l)
+ 4 \lam_1 \int^t_0 \int  | \nab ( \pd_r (\phi \frac{v_r}{r})) |^2  dy
+ 4 \lam_2 \int^t_0 \int   (  \pd_r (\phi \frac{v_r}{r}) )^2  dy\\
&\qquad + 4 \lam_1 \int^t_0 \int  | \nab ( \pd_3 (\phi \frac{v_r}{r})) |^2  dy
+ 4 \lam_2 \int^t_0 \int   ( \pd_3 (\phi \frac{v_r}{r}))^2  dy.
\eal
\ee Here the constant $C$ may have changed. 
  We need to bound the last 4 terms on the preceding inequality. For this purpose, 
we first need to prove the following localized version of a nice identity by Miao and 
Zheng. 
For any $q \in (1, \infty)$, there is a positive constant $c_q$ such that
\be
\lab{miao-zheng}
\al
&\Vert \nab ( \phi \pd_r \frac{v_r}{r}) \Vert_q \le c_q \Vert \Om \phi \Vert_q 
+ S(v_0, a, l),\\
&\Vert \nab^2 ( \phi \pd_r \frac{v_r}{r}) \Vert_q \le
 c_q \Vert \nab (\Om \phi) \Vert_q 
+ S(v_0, a, l).
\eal
\ee Here, as always $\Om=w_\theta/r$.
The proof of theses inequalities is given in

{\it step 3.3.}  From the identity 
\[
\Delta b = - \nab \times (w_\theta \-> {e_\theta}) = 
\left(\pd_3 ( w_\theta \frac{x_1}{r}),  \pd_3 (w_\theta \frac{x_2}{r}),  
\pd_1(w_\theta  \frac{x_1}{r}) - \pd_2 (w_\theta  \frac{x_2}{r})   \right),
\]and $b= v_r (\frac{x_1}{r}, \frac{x_2}{r}, 0) + v_3 (0, 0, 1)$, we see that
\be
\Delta (v_r \frac{x_1}{r}) = \pd_3 (x_1 \Om), \quad 
\Delta (v_r \frac{x_2}{r}) = \pd_3 (x_2 \Om).
\ee  Therefore
\be
\lab{ddvrx1phi}
\Delta (v_r \frac{x_1}{r} \phi) = \pd_3 (x_1 \Om \phi)- x_1 \Om \pd_3 \phi 
+ 2 \nab (v_r \frac{x_1}{r}) \nab \phi + v_r \frac{x_1}{r} \Delta \phi.
\ee Likewise
\be
\lab{ddvrx2phi}
\Delta (v_r \frac{x_2}{r} \phi) = \pd_3 (x_2 \Om \phi)- x_2 \Om \pd_3 \phi 
+ 2 \nab (v_r \frac{x_2}{r}) \nab \phi + v_r \frac{x_2}{r} \Delta \phi.
\ee  Inverting the Laplace operator, we infer
\be
\lab{vrx1phi}
v_r \frac{x_1}{r} \phi = \Delta^{-1} \pd_3 (x_1 \Om \phi) - \Delta^{-1}[ x_1 \Om \pd_3 \phi 
- 2 \nab (v_r \frac{x_1}{r}) \nab \phi - v_r \frac{x_1}{r} \Delta \phi],
\ee 
\be
\lab{vrx2phi}
v_r \frac{x_2}{r} \phi = \Delta^{-1} \pd_3 (x_2 \Om \phi) -
\Delta^{-1} [ x_2 \Om \pd_3 \phi 
- 2 \nab (v_r \frac{x_2}{r}) \nab \phi - v_r \frac{x_2}{r} \Delta \phi].
\ee Multiplying (\ref{vrx1phi}) by $x_1$, (\ref{vrx2phi}) by $x_2$ and taking the sum, 
we arrive at
\be
\lab{vrphi=1}
v_r \phi =  \Sigma^2_{i=1} \frac{x_i}{r} \Delta^{-1} \pd_3 (x_i \Om \phi)
- \Sigma^2_{i=1} \frac{x_i}{r} \Delta^{-1}[ x_i \Om \pd_3 \phi 
- 2 \nab (v_r \frac{x_i}{r}) \nab \phi - v_r \frac{x_i}{r} \Delta \phi].
\ee

Since $\phi$ is axially symmetric and $x_1/r = \cos \theta$, $x_2/r=\sin \theta$, 
we can write, for $i=1, 2$, that
\[
\nab (v_r \frac{x_i}{r}) \nab \phi = \frac{x_i}{r} (\pd_r v_r \pd_r \phi + 
\pd_3 v_r \pd_3 \phi ).
\] This turns (\ref{vrphi=1}) into
\be
\al
\lab{vrphi=2}
v_r \phi &= \Sigma^2_{i=1} \frac{x_i}{r} \Delta^{-1} \pd_3 (x_i \Om \phi)
- \Sigma^2_{i=1} \frac{x_i}{r} \Delta^{-1}( x_i f ),\\
&f \equiv \Om \pd_3 \phi 
- 2 \frac{\pd_r v_r}{r} \pd_r \phi
- 2 \frac{\pd_3 v_r}{r} \pd_3 \phi -  \frac{v_r}{r} \Delta \phi.
\eal
\ee Note the function $f$ is compactly supported, axially symmetric and point-wise bounded, due to 
the choice of the cut off function $\phi$.

According to \cite{MZ}, the following operator identity holds, at east when acting on compactly
supported functions, 
\be
\lab{mzdenshi1}
\Sigma^2_{i=1} \frac{x_i}{r} \Delta^{-1} x_i = r \Delta^{-1} - 2 \pd_r \Delta^{-2}.
\ee Since their proof is very sharp and cute, we repeat it here for completeness. 
Notice that
\[
\Sigma^2_{i=1} x_i [x_i, \Delta^{-1} ] = \Sigma^2_{i=1} x^2_i \Delta^{-1}
- \Sigma^2_{i=1} x_i \Delta^{-1} x_i= r^2 \Delta^{-1}
- \Sigma^2_{i=1} x_i \Delta^{-1} x_i.
\]Hence
\be
\lab{xixi1}
 \Sigma^2_{i=1} \frac{x_i}{r} \Delta^{-1} x_i
= r \Delta^{-1} - \Sigma^2_{i=1} \frac{x_i}{r} [x_i, \Delta^{-1} ].
\ee On the other hand 
\[
\Delta  [x_i, \Delta^{-1} ]  = \Delta (x_i \Delta^{-1}) - \Delta \Delta^{-1} x_i
= 2 \pd_i \Delta^{-1},
\]which implies
\[
[x_i, \Delta^{-1} ] =  2  \pd_i \Delta^{-2}.
\]Substituting this to the last term in (\ref{xixi1}), one obtains (\ref{mzdenshi1}).
Plugging  (\ref{mzdenshi1}) into the first identity in (\ref{vrphi=2}), we find that
\be
\lab{eqvr/r}
\frac{v_r}{r} \phi = (\Delta^{-1} \pd_3 - 2 \frac{\pd_r}{r} \Delta^{-2} 
\pd_3) (\Om \phi) -  (\Delta^{-1}  - 2 \frac{\pd_r}{r} \Delta^{-2} ) f.
\ee Recall that both $\Om \phi$ and $f$ are axially symmetric. When the operator 
$\frac{\pd_r}{r}$
acts on these functions, it can be written as 
\[
\frac{\pd_r}{r} =\Delta - \pd^2_r - \pd^2_3.
\]Plugging this into (\ref{eqvr/r}), we deduce 
\be
\lab{eqdvr/r}
\nab (\frac{v_r}{r} \phi ) = \Pi_1 (\Om \phi) +  \Pi_0 f,
\ee where $\Pi_1$ and $\nab \Pi_0$ are Riesz type singular integral operators
 that map $L^q$ to $L^q$, 
$q \in (1, \infty)$ and $\Pi_0$ is a smoothing integral operator. Since $f$ is bounded 
and compactly supported, this proves (\ref{miao-zheng}). We have used the fact that the 
gradient $\nab$ does not involve the derivative in $\->{e_\theta}$ direction, when acting on 
axially symmetric functions. 
\medskip

{\it step 3.4.}

Now we can take $q=2$ in  (\ref{miao-zheng}) and substitute it to (\ref{T1<4}) to 
obtain
\be
\lab{T1<8}
\al
 T_1
&\le \frac{1}{8} \int^t_0 \int |\nab J|^2 \phi^2 dy + C t S(v_0, a, l)
+ 4 \lam_1 c_2 \int^t_0 \int  | \nab (\Om \phi) |^2  dy
+ 4 \lam_2   c_2 \int^t_0 \int   (  \Om \phi )^2  dy\\
&\qquad + 4 \lam_1  c_2 \int^t_0 \int  | \nab (\Om \phi) |^2  dy
+ 4 \lam_2  c_2 \int^t_0 \int   ( \Om \phi)^2  dy.
\eal
\ee  This, together with (\ref{T2<}), yield
\be
\lab{T1T2<}
\al
T_1+ T_2 &\le (\frac{1}{8}+ 2 \lambda_1 + 9  \lam_1 c_2) 
\int^t_0 \int  (|\nab J|^2+ |\nab \Om |^2) \phi^2 dyds \\
&\qquad 
+ (\lambda_2 + 8 \lam_2  c_2) \int^t_0 \int [ (J \phi)^2 +  (\Om \phi)^2] dyds + C t S(v_0, a, l).
\eal
\ee In the above we have used the product formula $(\nab \Om ) \phi =
\nab (\Om \phi) - \Om \nab \phi$. This completes Step 3.

\medskip
{\it Step 4.  bounds on $T_3$.} 

Using $w_3 = \frac{1}{r} \pd_r ( r v_\theta)$, we compute 
\[
\al
&\int w_3 \pd_3 v_3 w_3 \phi^2 dy = \int \int^\infty_0 \pd_r ( r v_\theta) \pd_3 v_3 w_3 
\phi^2 dr dy_3 \\
&=-\int \int^\infty_0  r v_\theta \pd_r \pd_3 v_3 w_3 \phi^2 drdy_3 
- \int \int^\infty_0   r v_\theta  \pd_3 v_3 \pd_r w_3 \phi^2 drdy_3
- \int \int^\infty_0   r v_\theta  \pd_3 v_3 w_3 \pd_r \phi^2 drdy_3\\
&=-\int   v_\theta \pd_r \pd_3 v_3 w_3 \phi^2 dy
- \int    v_\theta  \pd_3 v_3 \pd_r w_3 \phi^2 dy
-  \int    v_\theta  \pd_3 v_3 w_3 \pd_r \phi^2 dy.
\eal
\] Next, using $w_r = -\pd_3 v_\theta$, we have
\[
\al
\int w_r \pd_r v_3 w_3 \phi^2 dy &= - \int \pd_3 v_\theta \pd_r v_3 w_3 \phi^2 dy\\
&=\int  v_\theta \pd_3 \pd_r v_3 w_3 \phi^2 dy +
\int  v_\theta \pd_r v_3 \pd_3 w_3 \phi^2 dy
+\int  v_\theta \pd_r v_3 w_3 \pd_3 \phi^2 dy.
\eal
\]Adding the previous two equalities and noting that the first terms on the 
right hand sides cancel, we obtain
\[
\al
T_3 &= -\int^t_0  \int    v_\theta  \pd_3 v_3 \pd_r w_3 \phi^2 dyds
- \int^t_0 \int    v_\theta  \pd_3 v_3 w_3 \pd_r \phi^2 dyds\\
&\qquad + \int^t_0 \int  v_\theta \pd_r v_3 \pd_3 w_3 \phi^2 dyds
+\int^t_0 \int  v_\theta \pd_r v_3 w_3 \pd_3 \phi^2 dyds.
\eal
\]As before, all terms involving derivatives of $\phi$ are bounded by $C t S(v_0, a, l)$.
Thus
\be
\lab{T3<1}
\al
T_3 &\le -\int^t_0  \int    v_\theta  \pd_3 v_3 \pd_r w_3 \phi^2 dyds
 + \int^t_0 \int  v_\theta \pd_r v_3 \pd_3 w_3 \phi^2 dyds
+C t S(v_0, a, l)\\
&\equiv I_1 + I_2 +C t S(v_0, a, l).
\eal
\ee

We will bound $I_1$ first.
By our condition on $v_\theta$,
\[
\al
I_1 &\le \frac{1}{8} \int^t_0  \int   |\pd_r w_3|^2 \phi^2 dyds +
      2 \int^t_0  \int    v^2_\theta  |\pd_3 v_3|^2 \phi^2 dyds\\
&\le \frac{1}{8} \int^t_0  \int   |\pd_r w_3|^2 \phi^2 dyds +
2 \lam_1 \int^t_0  \int      |\nab(\phi \pd_3 v_3)|^2  dyds
   +   2 \lam_2 \int^t_0  \int      |\pd_3 v_3|^2 \phi^2 dyds.
\eal
\]Consequently
\be
\lab{I1<1}
I_1 \le \frac{1}{8} \int^t_0  \int   |\pd_r w_3|^2 \phi^2 dyds +
3 \lam_1 \int^t_0  \int      |\nab \pd_3 v_3|^2 \phi^2  dyds
+ C t S(v_0, a, l, \lam_2).
\ee We need to bound the second term on the right hand side. To this end we call the
 relation for the full three dimensional velocity and vorticity: 
\[
- \Delta \pd_i v = \nab \times \pd_i w,
\]where $i=1, 2, 3$.  Using $\pd_i v \phi^2$ as a test function and integrate, we 
know that
\[
\al
&\int | \nab \pd_i v |^2 \phi^2 dy + \int \pd_j \pd_i v  \pd_i v \pd_j \phi^2  dy
=\int (\nab \times \pd_i w) \pd_i v \phi^2 dy\\
&=- \int (\nab \times  w) \pd_i \pd_i v \phi^2 dy 
- \int (\nab \times  w)  \pd_i v  \pd_i \phi^2 dy\\
&\le \frac{1}{2}\int | \nab \pd_i v |^2 \phi^2 dy + 
\frac{1}{2}\int |\nab \times  w|^2 \phi^2 dy  
- \int (\nab \times  w)  \pd_i v  \pd_i \phi^2 dy.
\eal
\] Since the terms involving derivatives of $\phi$ are bounded, this shows
\be
\lab{dd3v3}
\al
\int^t_0  \int     |\nab \pd_3 v_3|^2 \phi^2  dyds &\le 
\int^t_0\int | \nab \times w |^2 \phi^2 dyds + C t S(v_0, a, l)\\
&\le \int^t_0\int | \nab w |^2 \phi^2 dyds + C t S(v_0, a, l),
\eal
\ee and
\be
\lab{ddrv3}
\al
\int^t_0  \int      |\nab \pd_r v_3|^2 \phi^2  dyds &\le 
\int^t_0\int | \nab \times w |^2 \phi^2 dyds + C t S(v_0, a, l)\\
&\le \int^t_0\int | \nab w |^2 \phi^2 dyds + C t S(v_0, a, l).
\eal
\ee Here the constant $C$ may have changed when we drop the cross product,
which can be done through integration by parts  that produces extra bounded terms
involving $\nab \phi$.

Substituting (\ref{dd3v3}) into the second term on the right hand side of (\ref{I1<1}),
we reach 
\be
\lab{I1<2}
I_1 \le \frac{1}{8} \int^t_0  \int   |\pd_r w_3|^2 \phi^2 dyds +
3 \lam_1 \int^t_0  \int      |\nab w|^2 \phi^2  dyds
+ C t S(v_0, a, l, \lam_1, \lam_2).
\ee

Similarly, by our condition on $v_\theta$,
\[
\al
I_2 &\le \frac{1}{8} \int^t_0  \int   |\pd_3 w_3|^2 \phi^2 dyds +
      2 \int^t_0  \int    v^2_\theta  |\pd_r v_3|^2 \phi^2 dyds\\
&\le \frac{1}{8} \int^t_0  \int   |\pd_3 w_3|^2 \phi^2 dyds +
2 \lam_1 \int^t_0  \int      |\nab(\phi \pd_r v_3)|^2  dyds
   +   2 \lam_2 \int^t_0  \int      |\pd_r v_3|^2 \phi^2 dyds.
\eal
\]This with (\ref{ddrv3}) imply that
\be
\lab{I2<2}
I_2 \le \frac{1}{8} \int^t_0  \int   |\pd_3 w_3|^2 \phi^2 dyds +
3 \lam_1 \int^t_0  \int      |\nab w|^2 \phi^2  dyds
+ C t S(v_0, a, l, \lam_1, \lam_2).
\ee 

Substituting (\ref{I1<2}) and (\ref{I2<2}) into (\ref{T3<1}), we deduce the bound
for $T_3$, i.e.
\be
\lab{T3<2}
T_3 \le \frac{1}{8} \int^t_0  \int   |\nab w_3|^2 \phi^2 dyds +
6 \lam_1 \int^t_0  \int      |\nab w|^2 \phi^2  dyds
+ C t S(v_0, a, l, \lam_1, \lam_2).
\ee 
\medskip

{\it Step 5.  conclusion of the proof.} 

Combining (\ref{T1T2<}) with (\ref{T3<2}), we get
\be
\lab{T123<}
\al
T_1+ T_2 +T_3 &\le (\frac{1}{8}+ 2 \lambda_1 + 9  \lam_1 c_2) 
\int^t_0 \int  (|\nab J|^2+ |\nab \Om |^2) \phi^2 dyds \\
&\qquad 
+ (\lambda_2 + 8 \lam_2  c_2) \int^t_0 \int [ (J \phi)^2 +  (\Om \phi)^2] dyds +
\frac{1}{8} \int^t_0  \int   |\nab w_3|^2 \phi^2 dyds \\
&\qquad+
6 \lam_1 \int^t_0  \int      |\nab w|^2 \phi^2  dyds+ C t S(v_0, a, l, \lam_1, \lam_2).
\eal
\ee This, (\ref{R123<})  and (\ref{jowz}) together give
\[
\al
\int &\left(J^2+ \Om^2 +w^2_3 \right) \phi^2 dy \bigg |^t_0
+  
 \int^t_0 \int \left( | \nab J|^2 +|\nab \Om|^2 + | \nab w_3 |^2 \right) \phi^2 dyds\\
&\le  \frac{1}{4} \int^t_0 \int  \left( |\nabla J|^2  +
|\nabla \Om|^2 + |\nabla w_3|^2 \right) \phi^2 dyds \\ 
&\qquad  (\frac{1}{4}+ 4 \lambda_1 + 18  \lam_1 c_2) 
\int^t_0 \int  (|\nab J|^2+ |\nab \Om |^2) \phi^2 dyds \\
&\qquad 
+ 2(\lambda_2 + 8 \lam_2  c_2) \int^t_0 \int [ (J \phi)^2 +  (\Om \phi)^2] dyds +
\frac{1}{4} \int^t_0  \int   |\nab w_3|^2 \phi^2 dyds \\
&\qquad+
12 \lam_1 \int^t_0  \int      |\nab w|^2 \phi^2  dyds+ C t S(v_0, a, l, \lam_1, \lam_2).
\eal
\]Hence
\be
\lab{4+18}
\al
\int &\left(J^2+ \Om^2 +w^2_3 \right) \phi^2 dy \bigg |^t_0
+  \frac{1}{4}
 \int^t_0 \int \left( | \nab J|^2 +|\nab \Om|^2 + | \nab w_3 |^2 \right) \phi^2 dyds\\
&\le   ( 4 + 18  c_2) \lam_1
\int^t_0 \int  (|\nab J|^2+ |\nab \Om |^2) \phi^2 dyds +
 12 \lam_1 \int^t_0  \int      |\nab w|^2 \phi^2  dyds\\
&\qquad 
+ 2 \lam_2(1+ 8 c_2) \int^t_0 \int [ (J \phi)^2 +  (\Om \phi)^2] dyds
+ C t S(v_0, a, l, \lam_1, \lam_2).
\eal
\ee

There is still a little work to do, namely to bound the second term on the right hand side by the
left hand side.  Notice that $w$ is axially symmetric. Hence 
\[
\al
 |\nab w|^2 &= |\pd_r w_r|^2 + |\pd_r w_\theta|^2 + |\pd_3 w_r|^2 + |\pd_3 w_\theta|^2
+ |\nab w_3|^2\\
&= |\pd_r (J r)|^2 + |\pd_r (\Om r) |^2 +  r^2 |\pd_3 J|^2 + r^2 |\pd_3 \Om|^2
+ |\nab w_3|^2\\
&= |r \pd_r J + J|^2 + |r \pd_r \Om +\Om |^2 +  r^2 |\pd_3 J|^2 + r^2 |\pd_3 \Om|^2
+ |\nab w_3|^2\\
&\le 2 r^2 |\pd_r J|^2 + 2 J^2 + 2 r^2 |\pd_r \Om|^2
+ 2 \Om^2 +  r^2 |\pd_3 J|^2 + r^2 |\pd_3 \Om|^2
+ |\nab w_3|^2.
\eal
\]Hence
\[
|\nab w|^2 \le 2 r^2 (|\nab J|^2 + |\nab \Om|^2) + |\nab w_3|^2 + 2 (J^2+\Om^2).
\]Plugging this to the second term on the right hand side of (\ref{4+18}), we arrive at
\[
\al
\int &\left(J^2+ \Om^2 +w^2_3 \right) \phi^2 dy \bigg |^t_0
+  \frac{1}{4}
 \int^t_0 \int \left( | \nab J|^2 +|\nab \Om|^2 + | \nab w_3 |^2 \right) \phi^2 dyds\\
&\le   ( 28 + 18  c_2) \lam_1
\int^t_0 \int  (|\nab J|^2+ |\nab \Om |^2 +|\nab w_3|^2) \phi^2 dyds \\
&\qquad 
+ 2 [\lam_2(1+ 8 c_2)+24 \lam_1]  \int^t_0 \int [ (J \phi)^2 +  (\Om \phi)^2] dyds
+ C t S(v_0, a, l, \lam_1, \lam_2).
\eal
\]Here we have used the assumption that $r \le a \le 1$.
Choosing
\be
\lab{lam1}
\lam_1 = \frac{1}{4  ( 28 + 18  c_2)}.
\ee Here $c_2$ is given in (\ref{miao-zheng}) with $q=2$.
We reduce the last inequality to 
\[
\al
\int &\left(J^2+ \Om^2 +w^2_3 \right) \phi^2 dy \bigg |^t_0
\\
&\le 
 2 [\lam_2(1+ 8 c_2)+24 \lam_1]  \int^t_0 \int [ (J \phi)^2 +  (\Om \phi)^2] dyds
+ C t S(v_0, a, l, \lam_1, \lam_2).
\eal
\] By Gronwall's inequality 
\[
\int_{0 \le r \le a/2, \,  -l/2<y_3<l/2}
 \left((\frac{w_r}{r})^2+ (\frac{w_\theta}{r})^2 +w^2_3 \right) \phi^2(y, t) dy \le 
C(t, v_0, a, l, \lam_1, \lam_2).
\]By standard theory this is more than enough to imply the regularity of $v$ for all time.
The reason is that it implies $w$ is locally $L^{2, \infty}$ in any finite time. \qed
\medskip

Finally we verify the claim that $v_\theta$ is in the $\lam_1$ critical class for any fixed 
$\lam_1>0$,  if it satisfies 
$|v_\theta(x, t)| \le \frac{C}{r |\ln r|^{2+\e}}, \quad r<1/2$. 

Let $\psi=\psi(y, s)$ be any test function in Definition \ref{defcrit} with $a>0$ to be 
specified later. Fixing $s$, we compute
\[
\al
&\int \frac{\psi^2}{ r^2 |\ln r|^{2+\e}} dy 
= 2 \pi \int \int^\infty_0 \frac{1}{r |\ln r|^{2+\e}} \psi^2 dr dy_3\\
&=\frac{ 2 \pi}{1+\e} \int \int^\infty_0 \left( |\ln r|^{-1-\e} \right)' \psi^2 dr dy_3
=-\frac{ 2 \pi}{1+\e} \int \int^\infty_0 
 \frac{1}{|\ln r|^{1+(\e/2)}} \frac{2 \psi}{\sqrt{r}} \pd_r \psi  
\frac{1}{|\ln r|^{\e/2}} \sqrt{r} dr dy_3 \\
&\le \frac{ 2 \pi}{1+\e} \int \int^\infty_0 \frac{\psi^2}{r |\ln r|^{2+\e}} dr dy_3
+ \frac{ 2 \pi}{1+\e} \int \int^\infty_0 \frac{|\pd_r \psi|^2}{|\ln r|^{\e}} r dr dy_3\\
&\le \frac{1}{1+\e} \int  \frac{\psi^2}{r^2 |\ln r|^{2+\e}} dy
+ \frac{1}{1+\e} \int  \frac{|\pd_r \psi|^2}{|\ln r|^{\e}} dy.
\eal
\]Therefore
\[
\int \frac{\psi^2}{ r^2 |\ln r|^{2+\e}} dy \le \frac{1}{\e |\ln a|^\e} 
\int |\pd_r \psi|^2 dy,
\]which shows
\[
\int \left(\frac{|v_\theta|}{r} + v^2_\theta \right) \psi^2 dy \le \frac{C+C^2}{\e |\ln a|^\e} 
\int |\pd_r \psi|^2 dy.
\]Since $C$, $\e$ and $\lam_1$ are fixed positive numbers, we can always choose $a>0$ sufficiently 
small so that, for all $t \ge 0$, 
\[
\int^t_0\int \left(\frac{|v_\theta|}{r} + v^2_\theta \right) \psi^2 dyds \le \lam_1
\int^t_0\int |\pd_r \psi|^2 dyds.
\]Therefore $v_\theta$ is in the $\lam_1$ critical class.

\bigskip

{\it {\bf Acknowledgment}
 The author gratefully acknowledges the supports by
Siyuan Foundation through Nanjing University and by the Simons
Foundation.

He also wish to  thank Prof. Lei, Zhen and Mr.  Pan, Xinghong for discussions on the
problem.}


\begin{thebibliography}{99}

\bibitem{BZ:1} Burke Loftus, Jennifer and Zhang, Qi S. {\it A priori bounds for the vorticity of
axially symmetric solutions to the Navier-Stokes equations.} Adv.
Differential Equations 15 (2010), no. 5-6, 531-560.

\bibitem{CFZ:1} Hui Chen, Daoyuan Fang and Ting Zhang, {\it Regularity 
of 3D axisymmetric Navier-Stokes equations}, preprint 2015, to appear.

\bibitem{CKN:1}
 L. Caffarelli, R. Kohn, and L. Nirenberg,
 {\em Partial regularity of suitable weak solutions
 of the Navier-Stokes equations},
 Comm. Pure Appl. Math., 35 (1982),  771--831.

\bibitem{CL}
Dongho Chae and Jihoon Lee,
 {\em On the regularity of the axisymmetric
solutions of the Navier-Stokes equations},
 Math. Z.,  239  (2002),   645--671.

\bibitem{CSTY1}
Chiun-Chuan Chen, Robert M. Strain, Tai-Peng Tsai, and Horng-Tzer Yau,
{\em Lower bound on the blow-up rate of the axisymmetric
 Navier-Stokes equations},
Int. Math Res. Notices (2008), vol. 8, artical ID rnn016, 31 pp.

\bibitem{CSTY2} \bysame ,
{\em Lower bound on th blow-up rate of the axisymmetric
 Navier-Stokes equations II},Comm. P.D.E., 34(2009), no. 1--3, 203--232.

\bibitem{HL}
Thomas Y. Hou and Congming Li,
{\em Dynamic stability of the 3D axi-symmetric Navier-Stokes
equations with swirl},  Comm. Pure Appl. Math., 61 (2008)
661--697.

\bibitem{HLL}
Thomas Y. Hou, Zhen Lei, and Congming Li,
{\em Global reuglarity of the 3D axi-symmetric Navier-Stokes
equations with anisotropic data}, Comm. P.D.E., 33 (2008),
1622--1637.

\bibitem{JX}
 Quansen Jiu and Zhouping Xin,
 {\em Some regularity criteria on suitable weak
solutions of the 3-D incompressible axisymmetric
 Navier-Stokes equations}, Lectures on partial
differential equations,
 New Stud. Adv. Math., vol. \textbf{2},
Int. Press, Somerville, MA, 2003, pp. 119--139.


\bibitem{KNSS}
G. Koch, N. Nadirashvili, G. Seregin, and V. Sverak,
{\em Liouville theorems for the Navier-Stokes
equations and applications}, Acta Math. 203(2009), no. 1, 83--105.

\bibitem{L}
 O. A. Ladyzhenskaya,
{\em Unique global
solvability of the three-dimensional Cauchy
problem for the Navier-Stokes equations in the
 presence of axial symmetry}, Zap. Naucn. Sem.
Leningrad. Otdel. Math. Inst. Steklov. (LOMI)
\textbf{7} (1968), 155--177 (Russian).

\bibitem{LZ11} Z. Lei, and Q. S. Zhang, {\em A Liouville Theorem for the Axially-symmetric Navier-Stokes Equations}.
Journal of Functional Analysis, \textbf{261} (2011), 2323--2345.

\bibitem{LZ11-2} \bysame , {\em Structure of solutions of 3D Axi-symmetric Navier-Stokes Equations near Maximal Points}.
Pacific Journal of Mathematics, \textbf{254} (2011), no. 2, 335--344.

\bibitem{LMNP}
S. Leonardi, J. Malek, J. Necas, and M. Porkorny,
{\em On axially symmetric flows in $\reals^3
$}, Z. Anal. Anwendungen, \textbf{18} (1999),  639--649.

\bibitem{MZ} Miao, Changxing; Zheng, Xiaoxin, 
{\em On the global well-posedness for the Boussinesq system with horizontal dissipation.}
 Comm. Math. Phys. 321 (2013), no. 1, 33-67.

\bibitem{NP}
 Jiri Neustupa and Milan Pokorny,
 {\em An interior regularity criterion for an axially
 symmetric suitable weak solution to the Navier-Stokes equations},
 J. Math. Fluid Mech., 2 (2000),  381--399.

\bibitem{NP:2}
 Neustupa, J.; Pokorny, M. {\em
 Axisymmetric flow of Navier-Stokes fluid in the whole space with non-zero angular
 velocity component.}
 Proceedings of Partial Differential Equations and Applications (Olomouc, 1999).
 Math. Bohem. 126 (2001), no. 2, 469–481.

\bibitem{SS}
 G. Seregin and V. Sverak,
 {\em On type I singularities of
 the local axi-symmetric solutions of
the Navier-Stokes equations},
 Comm. P.D.E., 34(2009), no. 1--3, 171--201.

\bibitem{TX}
 Gang Tian and Zhouping Xin,
{\em One-point singular solutions to the Navier-Stokes
equations}, Topol. Methods Nonlinear Anal., 11
 (1998),
135--145.

\bibitem{UY}
 M. R. Ukhovskii and V. I. Yudovich,
{\em Axially symmetric flows of ideal and viscous fluids
filling the whole space}, J. Appl. Math. Mech., 32   (1968),
52--61.

\bibitem{Z}
 Qi S. Zhang,
 {\em A strong regularity result for parabolic equations},
 Comm. Math.
Phys., 244  (2004),   245--260.

\bibitem{ZZ} Zhang, Ping; Zhang, Ting, {\em
 Global axisymmetric solutions to three-dimensional Navier-Stokes system.}
 Int. Math. Res. Not. IMRN 2014, no. 3, 610-642. 


\end{thebibliography}
\end{document}